\documentclass[11pt]{amsart}
\usepackage{amsmath}
\usepackage{amssymb}
\usepackage{amsfonts}
\usepackage{graphicx}

\DeclareSymbolFont{AMSb}{U}{msb}{m}{n}
\DeclareMathSymbol{\N}{\mathbin}{AMSb}{"4E}
\DeclareMathSymbol{\Z}{\mathbin}{AMSb}{"5A}
\DeclareMathSymbol{\R}{\mathbin}{AMSb}{"52}
\DeclareMathSymbol{\Q}{\mathbin}{AMSb}{"51}
\DeclareMathSymbol{\I}{\mathbin}{AMSb}{"49}
\DeclareMathSymbol{\C}{\mathbin}{AMSb}{"43}

\theoremstyle{definition}

\theoremstyle{corollary}

\theoremstyle{example}

\theoremstyle{note}

\theoremstyle{notation}

\numberwithin{equation}{section}
\begin{document}
\title[Diagonal Sum of infinite image partition regular matrices]
{Diagonal Sum of infinite image partition regular matrices}
\author{Sourav Kanti Patra}
\address{Sourav Kanti Patra, Department of Mathematics, Ramakrishna Mission Vidyamandira,
Belur Math, Howrah-711202, West Bengal, India}
\email{souravkantipatra@gmail.com}
\author{Ananya shyamal}
\address{Ananya shyamal, Department of Electronics and Telecommunication, Engg., 
Army Institute of Technology, Alandi Road, Dighi Hills, Pune, Maharashtra 411015, India}
\email{shyamalananya18@gmail.com}


\begin {abstract}
A finite or infinite matrix $A$ is image partition regular
provided that whenever $\mathbb N$ is finitely colored, there must
be some $\vec{x}$ with entries from $\mathbb N$ such that all
entries of $A\vec{x}$ are in some color class. In [6], it was
proved that the diagonal sum of a finite and an infinite image partition
regular matrix is also image partition regular. It was also shown there
that centrally image partition regular matrices are closed
under diagonal sum. Using Theorem 3.3 of [2], one can conclude that diagonal
sum of two infinite image partition regular matrices may not be
image partition regular. In this paper we shall study the image partition
regularity of diagonal sum of some infinite image partition regular matrices. 
In many cases it will produce more infinite image partition regular matrices.

AMS subjclass [2010] : Primary : 05D10 Secondary : 22A15
\end{abstract}

\maketitle

\section {introduction}
Let us start this article with the following well known definition of 
image partition regularity.\\

\textbf{Definition1.1.} Let $u,v\in \mathbb{N}\cup\{\omega\}$ and let $A$
be a $u\times v$ matrix with entries from $\mathbb{Q}$. The matrix $A$
is image partition regular over $\mathbb{N}$ if whenever $r\in\mathbb{N}$ 
and $\mathbb{N}=\bigcup_{i=1} ^r C_i$, there exist $i\in\{1,2,.....,r\}$ 
and $\vec{x}\in\mathbb{N}^v$ such that $A\vec{x}\in C_{i}^u$.\\

Image partition regular matrices generalize many of the classical theorems of Ramsey
Theory. For example, Schur's Theorem [11] and the van der
Waerden's Theorem [12] are equivalent to say that the matrices
$\left(\begin{matrix} 1 & 0 \\ 0 & 1\\ 1 & 1
\end{matrix}\right)$ and for each $n\in \mathbb N$,
$\left(\begin{matrix} 1 & 0 \\ 1 & 1
\\\vdots & \vdots \\ 1 & n-1 \end{matrix}\right)$ are image
partition regular respectively.\\

It is well known that for finite matrices, image partition regularity 
behaves well with respect to central subsets of underlying semigroup.
Central sets were introduce by Furstenberg and
defined in terms of notion of topological dynamics. A nice
characterization of central sets in terms of algebraic structure
of $\beta\mathbb N$, the Stone-$\breve{C}$ech compactification of
$\mathbb N$ is given in Definition 1.3. Central sets are very
rich in combinatorial properties. The basic fact that we need
about central sets is given by the Central Sets Theorem, which is
due to
Furstenberg [3, Proposition 8.21] for the case $S = \mathbb Z$.\\

\textbf{Theorem 1.2.}(Central Sets Theorem)Let $(S,+)$ be a commutative semigroup.
Let $\tau$ be the set of sequences $\langle{y_t}\rangle_{t=1}^{\infty}$ in $S$.
Let $C$ be a subset of $S$ which is central and let $F\in\mathcal
P_{f}(\tau)$. Then there exist a sequence $\langle{a_t}\rangle_{t=1}^{\infty}$
in $S$ and a sequence $\langle{H_t}\rangle_{t=1}^{\infty}$ in $\mathcal
P_{f}(\mathbb N)$ such that for each $n\in\mathbb N$, $maxH_n <
minH_{n+1}$ and for each $L\in\mathcal P_{f}(\mathbb N)$ and each
$f\in F$, $\sum_{n\in L}(a_n + \sum_{t\in H_n} f(t))\in C$.\\

We shall present this characterization of central sets below, after
introducing the necessary background information.\\

Let $(S,\cdot)$ be an infinite discrete semigroup. Now the points
of $\beta S$ are taken to be the ultrafilters on $S$, the
principal ultrafilters being identified with the points of $S$.
Given $A\subseteq S$ let us set, $\bar{A}=\{p\in\beta S:A\in p\}$.
Then the set $\{\bar{A}:A\subseteq S\}$ will become a basis for a
topology on $\beta S$. The operation $\cdot$ on $S$ can be
extended to the Stone-$\breve{C}$ech compactification $\beta S$ of
$S$ so that $(\beta S,\cdot)$ is a compact right topological
semigroup (meaning that for any $p\in \beta S$, the function
$\rho_p:\beta S \rightarrow \beta S$ defined by $\rho_p(q)=q\cdot
p$ is continuous) with $S$ contained in its topological center
(meaning that for any $x\in S$, the function $\lambda_x:\beta S
\rightarrow \beta S$ defined by $\lambda_x(q)=x\cdot q$ is
continuous). Given $p,q\in\beta S$ and $A\subseteq S, A\in p\cdot
q$ if and only if $\{x\in
S:x^{-1}A\in q\}\in p$, where $x^{-1}A=\{y\in S:x\cdot y\in A\}$.\\

A nonempty subset $I$ of a semigroup $(T,\cdot)$ is called a left
ideal of $T$ if $T\cdot I\subseteq I$, a right ideal if
$I.T\subseteq I$, and a two-sided ideal (or simply an ideal) if it
is both a left and a right ideal. A minimal left ideal is a left
ideal that does not contain any proper left ideal. Similarly, we
can define minimal right ideal and smallest ideal. Any compact
Hausdorff right topological semigroup $(T,\cdot)$ has the unique
smallest two-sided ideal
$$\begin{array}{ccc}
K(T) & = & \bigcup\{L:L \text{ is a minimal left ideal of } T\} \\
& = & \,\,\,\,\,\bigcup\{R:R \text{ is a minimal right ideal of } T\}\\
\end{array}$$

Given a minimal left ideal $L$ and a minimal right ideal $R$ of
$T$, $L\cap R$ is a group, and in particular $K(T)$ contains an
idempotent. An idempotent that belongs to $K(T)$ is called a
minimal idempotent. We shall use the notation $\mathbb {N}^*$ for 
$\beta\mathbb N\setminus \mathbb N$\\

\textbf{Definition 1.3.} Let S be a semigroup and let $C\subseteq
S$. $C$ is called central in $S$ if there is some idempotent $p\in
K(\beta S)$ such that $C\in p$ (Definition 4.42, [7]).\\

In [10], R. Rado characterized kernel partition regular matrices 
in terms of computable condition called the column condition which 
generalized many of the classical theorem of Ramsey Theory.  
Like kernel partition regular matrices all the finite image
partition regular matrices can be described by a computable
condition called the first entries condition. In the following
theorem (Theorem 2.10, [5]) we see that central sets
characterize all finite image partition regular matrices.\\

\textbf{Theorem 1.4.} Let $u,v\in \mathbb N$ and let  $A$ be a $u
\times v$ matrix with entries from $\mathbb Q$. Then the following
statements are equivalent.\\
(a) $A$ is image partition regular.\\
(b) for every additively central subset $C$ of $\mathbb N$, there
exists $\vec{x}\in\mathbb{N}^v$ such that $A\vec{x}\in C^u$.\\

It is an immediate consequence of Theorem 1.4(b) that whenever $A$
and $B$ are finite image partition regular matrices, so is
$\left(\begin{matrix} A & 0 \\ 0 & B \end{matrix}\right)$, where
$0$ represents a matrix of appropriate size with zero entries.\\

Compared to finite image partition regular matrices, 
a little is known about the infinite one. In [4] and [6], the
notion of centrally image partition regular matrices were
introduced to extend the results of finite image partition regular
matrices to infinite image partition regular matrices. Now we recall 
the Definition 1.6(a) of [6].\\

\textbf{Definition 1.5.} Let $A$ be a $\omega\times\omega$ matrix
with entries from $\mathbb Q$. The matrix $A$ is centrally image
partition regular if for every central subset $C$ of $(\mathbb N,
+)$ there exists $\vec{x}\in\mathbb{N}^\omega$ such that
$A\vec{x}\in C^\omega$.\\

It follows immediately from the definition of centrally image partition 
regular matrices that the diagonal sum of any countable collection of
centrally image partition regular matrices is also image partition regular. 
This will produce more image partition regular matrices from the diagonal 
sum of centrally image partition regular matrices. 

We also recall the following definition(Definition 2.4, [6]).\\

\textbf{Definition 1.6.} Let $A$ be a finite or infinite matrix
with entries from $\mathbb{Q}$. Then $I(A)=\{p\in\beta\mathbb{N}$
: for every $P\in p$, there exists $\vec{x}$ with entries from
$\mathbb{N}$ such that all entries of $A\vec{x}$ are in $P\}$.\\

Following two lemmas study the algebraic structure of $I(A)$ 
with respect to $(\beta\mathbb{N},+)$ and $(\beta\mathbb{N},\cdot)$ 
respectively.\\

\textbf{Lemma 1.7.} Let $A$ be a matrix, finite or infinite with
entries from $\mathbb{Q}$.\\
(a) The set $I(A)$ is compact and $I(A)\neq\emptyset$ if and only
if $A$ is image partition regular.\\
(b) If $A$ is finite image partition regular matrix, then $I(A)$
is a sub-semigroup of $(\beta\mathbb{N},+)$.\\
\textbf{Proof.} See Lemma 2.5,[4].\\

\textbf{Lemma 1.8.} Let $A$ be a matrix, finite or infinite with
entries from $\mathbb{Q}$.\\
(a) If $A$ is an image partition regular matrix then $I(A)$ is a
left ideal of $(\beta\mathbb{N},\cdot)$.\\
(b) If $A$ is a finite image partition regular matrix then $I(A)$
is a two-sided ideal of $(\beta\mathbb{N},\cdot)$.\\
\textbf{Proof.} See Lemma 2.3, [9].\\

For a countable collection  $\langle{A_t}\rangle_{t=1}^{\infty}$ of 
matrices, let $B=\left(\begin{matrix} A_1 & 0 & 0 & \cdots \\ 
0 & A_2 & 0 & \cdots \\ 0 & 0 & A_3 & \cdots \\ \vdots & \vdots & \vdots & \cdots \\
\end{matrix}\right)$. Then $B$ is called the diagonal sum of the countable collection 
$\langle{A_t}\rangle_{t=1}^{\infty}$. In particular $\left(\begin{matrix} A_1 & 0 & 0 & \cdots & 0 \\ 
0 & A_2 & 0 & \cdots & 0 \\ 0 & 0 & A_3 & \cdots & 0 \\ 0 & 0 & 0 & \cdots & A_n\\
\end{matrix}\right)$ is the diagonal sum of finite collection of $n$ matrices 
$\langle{A_t}\rangle_{t=1}^{n}$. In this language we say, centrally image partition 
regular matrices are closed under diagonal sum to mean that diagonal sum of any two 
centrally image partition regular matrices is also centrally image partition regular.\\

In Lemma 2.3 of [6], diagonal sum of a finite and an infinite image partition regular matrix 
is also image partition regular. In section 2 and section 3 we shall introduce the notion of 
weak Milliken-Taylor system and substracted image partition regular matrices respectiviely to 
study their diagonal sum. At the end of this paper we shall see that the diagonal sum of a Weak
Milliken taylor matrix and a substracted image partition regular matrix is also image partition regular.

\section{Weak Milliken-Taylor system}

Milliken-Taylor system produces the class of infinite image partition 
regular matrices whose partition regularity does not come from centrally 
image partition regular matrices. For a compressed sequence (Definitio 2.1.(c))
$\vec{a}\in \mathbb{N}^m$, $m\in \mathbb{N}$, Milliken-Taylor study the
expression $a_1\cdot p+a_2\cdot p+...+a_m\cdot p$, where $p$ is be an additive 
idempotent of $\beta \mathbb N$ and $\vec{a}$=($a_1$, $a_2$, ... ,$a_m$). 
Being motivated by the Millike-Taylar system we have introduced the notion of 
weak Millike-Taylar system by considering the expression
$a_1\cdot p+a_2\cdot p+...+a_m\cdot p$ for any finite sequence,
$\vec{a}$=($a_1$, $a_2$, ... ,$a_m$)$\in \mathbb N$, $m\in\mathbb N$ and 
$p\in\beta\mathbb N$. Note that $p$ need not be an additive idempotent of 
$\beta\mathbb N$ in the expression $a_1\cdot p+a_2\cdot p+...+a_m\cdot p$ 
of weak Millike-Taylar system.\\

We now recall the Definition 2.1 of [8] and Definition 17.30 of [7] respectively.

\textbf{Definition 2.1.} Let $\vec{x}\in\omega^v$ where
$v\in\mathbb{N}\cup\{\mathbb N\}$.
Then\\
(a) $d(\vec{x})$ is the sequence obtained by deleting all
occurence of $0$ from $\vec{x}$.\\
(b) $c(\vec{x})$ is the sequence obtained by deleting every digit
in $d(\vec{x})$ which is equal to its predecessors and\\
(c) $\vec{x}$ is a compressed sequence if and only if
$\vec{x}=c(\vec{x})$.\\

\textbf{Definition 2.2.} Let $\vec{a}\in \mathbb{N}^m$ be a compressed sesuence, $m\in \mathbb{N}$ 
and a sequence $\langle{x_t}\rangle_{t=1}^{\infty}$ in $\mathbb{N}$,
$MT(\vec{a}, \langle{x_t}\rangle_{t=1}^{\infty})=\{\sum_{i=1}^{m} a_i\sum_{t\in F_i}x_{t}:F_1, F_2,...,F_m\in
\mathcal{P}_f(\mathbb{N})$ and $F_1<F_2<...<F_m\}$. Where $F<G$ means $\max F<\min G$ for $F,G\in \mathcal{P}_f(\mathbb{N}$)\\

Similarly we can define Product Milliken-taylor system as follows.

\textbf{Definition 2.3.} Let $\vec{a}\in \mathbb{N}^m$ be a finite sequence, $m\in \mathbb{N}$ 
and a sequence $\langle{x_t}\rangle_{t=1}^{\infty}$ in $\mathbb{N}$,
$PMT(\vec{a}, \langle{x_t}\rangle_{t=1}^{\infty})=\{\sum_{i=1}^{m} a_i\prod_{t\in F_i}x_{t}:F_1, F_2,...,F_m\in
\mathcal{P}_f(\mathbb{N})$ and $F_1<F_2<...<F_m\}$. Where $F<G$ means $\max F<\min G$ for $F,G\in \mathcal{P}_f(\mathbb{N}$)\\

we also recall the Definition 5.13 (b) of [7].

\textbf{Definition 2.4.} Let $\langle{x_t}\rangle_{t=1}^{\infty}$ be a sequence
in $\mathbb{N}$. A sequence $\langle{y_t}\rangle_{t=1}^{\infty}$ in $\mathbb{N}$
is said to be a sum-subsystem of $\langle{x_t}\rangle_{t=1}^{\infty}$ if there
exists a sequence $\langle{H_t}\rangle_{t=1}^{\infty}$ of finite subsets of
$\mathbb{N}$ with $\max H_t<\min H_{t+1}$ for all $t\in\omega$
such that $y_t=\sum_{s\in H_t}x_s$.\\

We shall now define Weak Milliken-taylor System.

\textbf{Definition 2.5} Let $\vec{a}\in \mathbb{N}^m$, $m\in \mathbb{N}$ 
and a sequence $\langle{x_t}\rangle_{t=1}^{\infty}$ in $\mathbb{N}$,
$WMT(\vec{a}, \langle{x_t}\rangle_{t=1}^{n})=\{\sum_{i=1}^{m} a_ix_{t_i}:t_1<t_2<...<t_m\}$\\

\textbf{Definition 2.6.} A $\omega\times\omega$ matrix $M$ with
entries from $\omega$ is said to be a Milliken-Taylor matrix if\\
(a) each row of $M$ has only finitely many nonzero entries and\\
(b) there exists a finite compressed sequence $\vec{a}\in\omega^m,
m\in\mathbb{N}$ such that $M$ consists of all possible  row
vectors $\vec{r}\in\omega^{\omega}$ for which $c(\vec{r})=\vec{a}$.\\

Similarly we now define Weak Milliken-Taylor matrix.

\textbf{Definition 2.7.} A $\omega\times\omega$ matrix $M$ with
entries from $\omega$ is said to be a weak Milliken-Taylor matrix if\\
(a) each row of $M$ has only finitely many non-zero entries and\\
(b) there exists a finite sequence $\vec{a}\in\omega^m,
m\in\mathbb{N}$ such that $M$ consists of all possible row
vectors $\vec{r}\in\omega^{\omega}$ for which $d(\vec{r})=\vec{a}$.\\

\textbf{Theorem 2.8.} Weak Milliken-Taylor matrices are centrally image partition 
regular.\\
\textbf{Proof.} Let $M$ be a weak Milliken-Taylor matrix. Then by Definition 2.7, 
there exists a finite sequence $\vec{a}\in\mathbb{N}^m,
m\in\mathbb{N}$ such that $M$ has all possible row
vectors $\vec{r}\in\omega^{\omega}$ for which $d(\vec{r})=\vec{a}$. 
Let $\sum_{i=1}^{m}a_i=k$ where $\vec{a}$=($a_1$ $a_2$ .... $a_m$). 
Given a central set $C$, simply pick $d\in\mathbb N$ such that $dm\in C$, 
which one can do because for each $n\in\mathbb N$, $n\mathbb N$ is a
member of every idempotent by Lemma 6.6 of [5]. Then let $x_i = d$ for 
each $i\in\omega$. Let $\vec{x}=\left(\begin{matrix} x_0 \\ x_1 \\ x_2 \\
\vdots\end{matrix}\right)$. Then $A\vec{x}\in C^\omega$.\\

\textbf{Remark 2.9.} As a consequence of the above theorem we can conclude that 
the diagonal sun of any sequence of weak Milliken-Taylor matrices is also 
centrally image partition regular, in particular they are image  partition regular.\\

\textbf{Lemma 2.10.} Let $p \in \mathbb {N}^*$, let$A \in p$ and $\langle{a_t}\rangle_{t=1}^{m}$ be any
finite sequence in $\mathbb{N}$ and let $B \in a_1\cdot p+a_2\cdot p+ ... +a_m\cdot p$. There is a one
to one sequence $\langle{x_t}\rangle_{t=1}^{\infty}\in A $ such that WMT($\vec{a},
\langle{x_t}\rangle_{t=1}^{\infty}) \subset B$.\\
\textbf{Proof} Assume first that m=1. Then $a_1^{-1}B \in p$,  so $A\cap a_1^{-1}B \in p$. 
Now choose an infinite one to one sequence $\langle{x_t}\rangle_{t=1}^{\infty}$ in $A\cap a_1^{-1}B $.
This says precisely that WMT($\vec{a}, \langle{x_t}\rangle_{t=1}^{\infty}) \subset B$ for some one to one 
sequence $\langle{x_t}\rangle_{t=1}^{\infty}$ in A. Assume now that $ m\geq2$ and notice that\\
$\{x \in \mathbb{N}:-x+B \in a_1\cdot p+a_2\cdot p+...+a_m\cdot p\} \in a_1\cdot p$\\
so that $\{x \in \mathbb{N}:-a_1x+B \in a_1\cdot p+a_2\cdot p+...+a_m\cdot p\}\in p$.\\
Let $B_1=\{x \in \mathbb{N}:-a_1x+B \in a_1\cdot p+a_2\cdot p+...+a_m\cdot p\}\cap A$. Now pick $x_1 \in B_1$.
Inductively, let $n \in \mathbb{N}$ and assume that we have choosen one to one sequence
$\langle{x_t}\rangle_{t=1}^{n}$ in $\mathbb{N}$, $\langle{B_t}\rangle_{t=1}^{n}$ in $p$, 
so that for each $ r \in \{1, 2, ..., n\}$\\
(I)If $t \in \{1, 2, ..., n\}$, $ x_t \in B_t$.\\
(II)If $r<n$, then $ B_{r+1} \subset B_r$.\\
(III)If $l\in \{1,2,...,m-1\}$, $t_1, t_2, ..., t_l \in\{1,2,...,r\}$ and $t_1<t_2<.....<t_l$, 
then $-\sum_{i=1}^{l} a_ix_{t_i}+B \in a_{l+1}.p+a_{l+2}.p+.....+a_m.p$\\
(IV)If $t_1, t_2,..., t_{m-1} \in\{1, 2, ..., r\}, t_1<t_2<.....<t_{m-1}$ and $r<n$, 
then $B_{r+1}\subset {a_m}^{-1}(-\sum_{i=1}^{l} a_ix_{t_i}+B)$.\\
(V)If $ l\in \{1,2,...,m-2\}, t_1, t_2, ..., t_l \in\{1,2,...,r\}, t_1<t_2<.....<t_l$ and $r <n$, 
then $ B_{r+1}\subset \{x \in \mathbb{N}:-a_{l+1}x +(-\sum_{i=1}^{l} a_ix_{t_i}+B) \in a_{l+2}\cdot p+a_2\cdot p+...+a_m\cdot p\}$.\\ 
At n=1, hypothesis (I) holds directly, hypothesis (II),(IV) and (V) are vacuous, 
and hypothesis (III) says that $-a_1x_1 +A \in a_2.p+a_3.p+...+a_m.p\}$  which is 
true because $x_1 \in B_1$.\\
For $ l\in \{1,2,...,m-1\}$,
let\\
$\mathcal{F}_l=\{(t_1,t_2,...,t_l) :t_1,t_2,...,t_l \in\{1,2,...,n\}$
and $t_1<t_2<.....<t_l\}$.\\
If $(t_1,t_2,...,t_l) \in \mathcal{F}_{m-1}$, then by hypothesis (III) we have
$-\sum_{i=1}^{l} a_ix_{t_i}+B \in a_m\cdot p$, so that 
${a_m}^{-1}(-\sum_{i=1}^{l} a_ix_{t_i}+B)\in p$. If $l\in \{1, 2, ..., m-2\}$ 
and $(t_1, t_2, ..., t_l) \in \mathcal{F}_l$, 
we have by (III) that $-\sum_{i=1}^{l} a_ix_{t_i}+B \in \{a_{l+1}\cdot p+a_2\cdot p+...+a_m\cdot p\}$, 
so then $ \{x \in \mathbb{N}:-a_{l+1}x +(-\sum_{i=1}^{l} a_ix_{t_i}+B) \in a_{l+2}\cdot p+......+a_m\cdot p\}\in p$\\
Now let $B_{n+1}=B_n \bigcap_{k=1}^{n}(-x_k+B_k)\bigcap\limits_ {\{(t_1,t_2,...,t_{m-1}) \in \mathcal{F}_{m-1}\}} {a_m}^{-1}(-\sum_{i=1}^{l} a_ix_{t_i}+B)\bigcap_{l=1}^{m-2}\bigcap \limits_ {\{(t_1,t_2,...,t_l) \in \mathcal{F}_l\}}\{x \in \mathbb{N}:-a_{l+1}x +(-\sum_{i=1}^{l} a_ix_{t_i}+B) \in
a_{l+2}.p+a_2.p+...+a_m.p\}$\\ Then we have $ B_{n+1} \in p$. Since $p\in \mathbb{N}^*, B_{n+1}$
is an infinite set. Now choose $x_{n+1}\in B_{n+1}\setminus \{x_1,x_2,...,x_n\}$ . Hypothesis (I) 
and (II) holds trivially. Hypothesis (IV) and (V) holds directly.To verify hypothesis (III) 
let $l\in \{1, 2, ..., m-1\}$ and $t_1<t_2<...<t_l \in \{1, 2, ..., n+1\}$. If $l=1$, by hypothesis (I),(II) 
we have $x_{t_1}\in B_1$, so then $-a_1x_{t_1} +B \in a_2\cdot p+......+a_m\cdot p$ as required. So assume $l>1$, 
then $x_{t_l} \in B_{t_l}\subset B_{t_l+1} \subset \{x \in \mathbb{N}:-a_lx +(-\sum_{i=1}^{l-1} a_ix_{t_i}+B) \in a_{l+1}\cdot p+a_{l+2}\cdot p+...+a_m\cdot p\}$\\
(by hypothesis(V) for $r=t_l$). So $-\sum_{i=1}^{l} a_ix_{ti}+A \in a_{l+1}.p+a_2.p+...+a_m.p$ as required.
The induction being complete, let $t_1<t_2<...<t_m \in \mathbb{N}$. 
Then $x_{t_m} \in B_{t_m}\subset B_{t_m+1} \subset {a_m}^{-1}(-\sum_{i=1}^{m-1} a_ix_{t_i}+B)$\\  So $\sum_{i=1}^{l} a_ix_{t_i} \in B$ 
as required.\\

\textbf{Theorem 2.11.} Let $\vec{a}=\langle a_1, a_2, ..., a_m\rangle $ be a finite sequence
in  $\mathbb{N}$.Let $p.p=p\in\bigcap_{k=1}^\infty\overline{FP(\{x_t\}_{t=k}^\infty)}$, 
and let  $A \in a_1.p+a_2.p+...+a_m.p$. There is a product subsystem $\langle x_t\rangle _{t=1}^\infty$
of $\langle y_t\rangle _{t=1}^\infty$  such that PMT( $\vec{a},\langle x_t\rangle _{t=1}^\infty \rangle \subset A$\\
\textbf{proof} Imitate the proof of theorem 17.31 of [7].\\
There is a partial converse to the Therem 2.11 bellow. 
Notice that $p$ is not required to be multiplicative idempotent in the following theorem.\\

\textbf{Theorem 2.12.} Let $\vec{a}=\langle a_1, a_2, ..., a_m\rangle$ be a finite sequence in  $\mathbb{N}$  
and let  $p\in\bigcap_{k=1}^\infty\overline{FP(\{x_t\}_{t=k}^\infty)}$. Then 
PMT( $\vec{a},\langle x_t\rangle _{t=k}^\infty)  \in a_1.p+a_2.p+...+a_m.p$\\
\textbf{proof}Imitate the proof of theorem 17.32 of [7].\\

\textbf{Theorem 2.13.} Let $a_0, a_1, ..., a_n,b_0, b_1, ..., b_m\in \mathbb{N}$ be so that for any 
$i \in\{1,2,...,n-1\}$ and $j\in\{1,2,...,m-1\}$   $ a_i \neq a_{i+1} ,b_j \neq b_{j+1}$. 
Let $p$ and $q$ are idempotents in $(\beta\mathbb{N},+)$. Suppose that 
$a_0+a_1\cdot p+ a_2\cdot p+a_3\cdot p+...+a_n\cdot p=b_0+b_1\cdot q+ b_2\cdot q+b_3\cdot q+...+b_m\cdot q$, then $a_0=b_0 $ and
for any $i\in\{1,2,...,n\},a_i=\lambda b_i$ for some $\lambda \in \mathbb{Q}\setminus \{0\}$.\\
\textbf{proof} If possible let $ a_0 \neq b_0$. Choose $d \in \mathbb{N}$ such that
$d \nmid  a_0 - b_0$. Let $\gamma :\mathbb{N}\rightarrow \mathbb{Z}/n\mathbb{Z}$ denote
the canonical map. Observe that $\gamma :(\mathbb{N},+)\rightarrow (\mathbb{Z}/n\mathbb{Z},+)$
and $\gamma :(\mathbb{N},\cdot)\rightarrow (\mathbb{Z}/n\mathbb{Z},\cdot)$ are semigroup homomorphisms. 
Also let $\tilde{\gamma}:\beta\mathbb{N} \rightarrow \mathbb{Z}/d\mathbb{Z}$ be the continuous
extension of $\gamma$. Then $\tilde{\gamma}:(\beta\mathbb{N},+) \rightarrow (\mathbb{Z}/d\mathbb{Z},+)$
and $\tilde{\gamma}:(\beta\mathbb{N},\cdot) \rightarrow (\mathbb{Z}/d\mathbb{Z},\cdot)$ both are semigroup
homomorphisms by corollary 4.22,[7]. Thus $\tilde{\gamma}(p)= \tilde{\gamma}(p) + \tilde{\gamma}(p)$
and so $\tilde{\gamma}(p)=0$. Similarly $\tilde{\gamma}(q)=0$. 
Now $\tilde{\gamma}(a_0+a_1\cdot p+a_2\cdot p+a_3\cdot p+...+a_n\cdot p)=\tilde{\gamma}(b_0+b_1\cdot q+ b_2\cdot q+b_3\cdot q+...+b_m\cdot q)$. 
Therefore $\tilde{\gamma}(a_0)+ (\sum_{i=1}^{n} \tilde{\gamma}(a_i)).\tilde{\gamma}(p)=\tilde{\gamma}(b_0)+ (\sum_{i=1}^{m} \tilde{\gamma}(b_i)).\tilde{\gamma}(q)$ .\\
Hence $ \gamma(a)=\gamma(b)$  $[\therefore    \tilde{\gamma}(p)  = \tilde{\gamma}(q)=0 $ and $ \tilde{\gamma}|_{\mathbb{N}}=\gamma]$. 
Hence $\bar{a}=\bar{b}$ i.e.$ d\mid (a_0 -b_0)$, which is a contradiction to  our assumption $d \nmid  a_0 - b_0$. 
Thus $a_0=b_0$.Hence $a_1.p+ a_2.p+a_3.p+...+a_n.p=b_1.q+ b_2.q+b_3.q+...+b_m.q$. 
So by Theorem 3.3, [2], we have$m=n$ and  $a_i=\lambda b_i$.\\

As a consequence of Theorem 17.27, [7] we can say that there do not exist multiplicative idempotents $p$ and $q$ of 
$(\beta\mathbb{N},\cdot)$ such that $p=q+q+...+q$($m$ times, $m\in \mathbb N-{1}$). This allows us to raise the following 
question.\\

\textbf{Question 2.14.} Let $a_1\cdot p+ a_2\cdot p+a_3\cdot p+...+a_n\cdot p=+b_1\cdot q+ b_2\cdot q+b_3\cdot q+...+b_m\cdot q$ 
for some multiplicative idempotents p,q in $(\beta\mathbb{N},\cdot)$. Can we say that $m=n$, $a_i=b_i$ for $i=1, 2, ..., m$ and 
$p=q$ \\

\section{substracted image partition regularity}
Let us start this section with the definition of substracted image partition regular matrices. 

\textbf{Definition 3.1.} Consider the following conditions for a
$\omega\times\omega$ matrix $A$.

(1) no row of $A$ is $\vec{0}$.

(2) for each $i\in\omega$, $\{j\in\omega:a_{ij}\neq 0\}$ is
finite.

(3) If $\vec{c_0},\vec{c_1},\vec{c_2},.....$ be the columns of
$A$, there exist $n\in\omega$, $k\in \mathbb N$ such that all the rows of
($\vec{c_{n}}$ $\vec{c_{n+1}}$
$\vec{c_{n+2}}$.....$\vec{c_{n+k-1}}$) are precisely the rows of a
finite image partition regular matrix and the remaining columns
form a  image partition regular matrix.

(4) If $\vec{c_0},\vec{c_1},\vec{c_2},.....$ be the columns of
$A$, there exist $n\in \omega$, $k\in \mathbb{N}$ such that all the rows of
($\vec{c_{n}}$ $\vec{c_{n+1}}$
$\vec{c_{n+2}}$.....$\vec{c_{n+k-1}}$) are precisely the rows of a
finite image partition regular matrix and the remaining columns of
$A$ form a centrally image partition regular matrix.\\

We call a matrix $A$ to be 'substructed image partition
regular' if it satisfies (1) (2) and (3) and also we call $A$ to
be 'substructed centrally image partition regular' if it satisfies
(1) (2) and (4). \\

Note that both of these matrices defined in Definition 3.1 are
image partition regular. Also observe that substracted centrally image
partition regular matrix is centrally image partition regular.\\

\textbf{Definition 3.2.} Let $M$ be an infinite image partition 
regular matrix with entries from $\mathbb Q$. A
$\omega\times\omega$ matrix $A$ with entries from $\mathbb Q$ 
is said to be $M$-substracted image partition regular matrix if

(1) no row of $A$ is $\vec{0}$.

(2) for each $i\in\omega$, $\{j\in\omega:a_{ij}\neq 0\}$ is
finite.

(3) If $\vec{c_0},\vec{c_1},\vec{c_2},.....$ be the columns of
$A$, there exist $n\in \omega$, $k\in \mathbb{N}$ such that all
the rows of ($\vec{c_{n}}$ $\vec{c_{n+1}}$
$\vec{c_{n+2}}$.....$\vec{c_{n+k-1}}$) are precisely the rows of a
finite image partition regular matrix and the remaining columns of
$A$ are precisely the columns of $M$.\\

Note that this is a particular class of substracted image partition 
regular matrices. As a consequence of the following theorem we can
conclude that the diagonal sum of a sequence of $M$-substracted 
image partition regular matrices is also  image partition regular.\\

\textbf{Theorem 3.3.} Let $M$ be an infinite image partition regular 
matrix and $\mathcal{M}$ be a collection of $M$-substracted image partition
regular matrices. Then $\bigcap_{A\in\mathcal{M}}I(A)\neq\emptyset$.\\
\textbf{Proof.} Let $\mathcal{F}$ be the set of all
finite image partition regular matrices with entries from $\mathbb Q$. 
Then by Lemma 1.7(b), for each $B\in\mathcal{F}$, $I(B)$ is a
two-sided ideal of $(\beta\mathbb N,\cdot)$ and so is
$\bigcap_{B\in\mathcal{F}}I(B)$. Also by Lemma 1.7(a) $I(M)$ is a 
left ideal of $(\beta\mathbb N,\cdot)$. Thus
$I(M)\cap(\bigcap_{B\in\mathcal{F}}I(B))\neq\emptyset$. Choose 
$p\in I(M)\cap(\bigcap_{B\in\mathcal{F}}I(B))$. Now let $A\in \mathcal{M}$.
If $\vec{c_0},\vec{c_1},\vec{c_2},.....$ be the
columns of $A$ then choose $n\in\omega$, $k\in \mathbb N$ and take
$F$=($\vec{c_{n}}$ $\vec{c_{n+1}}$.....$\vec{c_{n+k-1}}$) and
$M$ be the remaining columns of $A$ as in the definition of
$M$-substracted image partition regular matrix. Take
$C$=($F$ $M$) and observe that $I(A)=I(C)$. Now all the rows
of $F$ is precisely the rows of a finite image partition regular
matrix (i.e $F$ has infinitely many repeated rows precisely
coming from a particular finite image partition regular matrix.). 
Let $V\in p+p$.
Then $\{x\in\mathbb{N} : -x+V\in r\}\in q$. Now since $q\in
I(F)$, there exists $\vec{x^{(1)}}\in \mathbb{N}^k$ such that
$y_i\in\{x\in\mathbb{N} : -x+V\in r\}$ for all $i\in \omega$ where
$\vec{y}=A_1\vec{x^{(1)}}$ and
$\vec{y}=\left(\begin{matrix} y_0 \\ y_1 \\ y_2 \\
\vdots\end{matrix}\right)$. Hence $-y_i+V\in r$ for all
$i\in\omega$. Also observe that $\{y_i : i\in \omega\}$ is finite.
Thus $\bigcap_{i\in\omega}(-y_i+V)\in r$. Since $p\in I(M)$ there exists $\vec{x^{(2)}}\in
\mathbb{N}^{\omega}$ such that whenever $\vec{z}=A_2\vec{x^{(2)}}$
we have $z_j\in\bigcap_{i\in\omega}(-y_i+V)$ for all $j\in\omega$
where $\vec{z}=\left(\begin{matrix} z_0 \\ z_1 \\ z_2 \\
\vdots\end{matrix}\right)$. So $y_i+z_j\in V$ for all
$i,j\in\omega$. Now let $\vec{x}= \left(\begin{matrix}
\vec{x^{(1)}}\\ \vec{x^{(2)}} \end{matrix}\right)$. Then
$C\vec{x}=A_1\vec{x^{(1)}}+A_2\vec{x^{(2)}}=\vec{y}+\vec{z}$.
Therefore $C\vec{x}\in V^\omega$ and hence $p+p\in I(C)=I(A)$ 
Therefore $p+p\in I(A)$ for 
all $A\in \mathcal{M}$. Hence $p+p\in \bigcap_{A\in\mathcal{M}}I(A)$ and so
$\bigcap_{A\in\mathcal{M}}I(A)\neq\emptyset$.\\

In Theorem 2.16 of [9], it was shown that the diagonal sun of a substracted 
centrally image partition regular matrix and a Milliken-Taylor matrix 
is also image partition partion regular. In the following we are proposed 
to show that that the diagonal sun of a substracted  image partition regular
matrix and a weak Milliken-Taylor matrix is also image partition partion regula.\\

\textbf{Theorem 3.4.} Let  $A$ be a subtracted image partition regular 
matrix and $B$ be a Weak Milliken-Taylor matrix. Then 
$\left(\begin{matrix} A & 0 \\ 0 & B \end{matrix}\right)$ is image
partition regular.\\
\textbf{Proof} If $\vec{c_0},\vec{c_1},\vec{c_2},.....$ be the
columns of $A$ then choose $n\in\omega$, $k\in \mathbb N$ and take
$A_1$=($\vec{c_{n}}$ $\vec{c_{n+1}}$.....$\vec{c_{n+k-1}}$) and
$A_2$ to be the remaining columns of $A$ as in the definition of
subtracted image partition regular matrix. Take
$C$=($A_1$ $A_2$) and observe that $I(A)=I(C)$. Now all the rows
of $A_1$ is precisely the rows of a finite image partition regular
matrix (i.e $A_1$ has infinitely many repeated rows precisely
coming from a particular finite image partition regular matrix.).
Therefore $I(A_1)$ is a sub-semigroup of $(\beta\mathbb N,+)$ and
is a two-sided ideal of $(\beta\mathbb N,\cdot)$ by Lemma 1.7(b) 
and Lemma 1.8(b) respectively.
Since $B$ is a Weak Milliken-Taylor matrix, there is a finite 
sequence $\langle a_t\rangle _{t=1}^m$ $m\in\mathbb N$ 
such that $c(\vec{r})=\vec{a}$ for each row $\vec{r}$ of $M$.
Also $I(A_2)$ is a left ideal of $(\beta\mathbb N,\cdot)$. Now 
chooose $p\in I(A_1)\cap I(A_2)$. Take
$q = a_1\cdot p+a_2\cdot p+.....+a_{m-1}\cdot p$ and $r=a_m\cdot p$. 
Then $q\in I(A_1)$ and $r\in I(A_2)$. Let $V\in q+r$.
Then $\{x\in\mathbb{N} : -x+V\in r\}\in q$. Now since $q\in
I(A_1)$, there exists $\vec{x^{(1)}}\in \mathbb{N}^k$ such that
$y_i\in\{x\in\mathbb{N} : -x+V\in r\}$ for all $i\in \omega$ where
$\vec{y}=A_1\vec{x^{(1)}}$ and
$\vec{y}=\left(\begin{matrix} y_0 \\ y_1 \\ y_2 \\
\vdots\end{matrix}\right)$. Hence $-y_i+V\in r$ for all
$i\in\omega$. Also observe that $\{y_i : i\in \omega\}$ is finite.
Thus $\bigcap_{i\in\omega}(-y_i+V)\in r$. Since $r\in I(A_2)$ there exists $\vec{x^{(2)}}\in
\mathbb{N}^{\omega}$ such that whenever $\vec{z}=A_2\vec{x^{(2)}}$
we have $z_j\in\bigcap_{i\in\omega}(-y_i+V)$ for all $j\in\omega$
where $\vec{z}=\left(\begin{matrix} z_0 \\ z_1 \\ z_2 \\
\vdots\end{matrix}\right)$. So $y_i+z_j\in V$ for all
$i,j\in\omega$. Now let $\vec{x}= \left(\begin{matrix}
\vec{x^{(1)}}\\ \vec{x^{(2)}} \end{matrix}\right)$. Then
$C\vec{x}=A_1\vec{x^{(1)}}+A_2\vec{x^{(2)}}=\vec{y}+\vec{z}$.
Therefore $C\vec{x}\in V^\omega$ and hence $q+r\in I(C)=I(A)$. 
Also since $p\in \mathbb N^*$ and $V\in a_1\cot p+a_2\cdot p+...+a_m\cot p$, 
by Lemma 2.10 there is a one-to-one
sequence $\langle x_t\rangle _{t=1}^\infty$ in $\mathbb N$ such that
$WMT(\vec{a}, \langle x_t\rangle _{t=1}^\infty)\subseteq V$. Let
$\vec{x}=\left(\begin{matrix} x_0 \\ x_1 \\ x_2 \\
\vdots\end{matrix}\right)$. Then $B\vec{x}\in V^\omega$. Thus
$q+r\in I(B)$. Therefore $q+r\in I(A)\bigcap I(B)$. Since
$I(A)\bigcap I(B)\neq\emptyset$, it follows immediately from the
definition of $I(A)$ and $I(B)$ that $\left(\begin{matrix} A & 0 \\
0 & B \end{matrix}\right)$ is image partition regular.\\

\textbf {Acknowledgement.}  The authors are grateful to Prof.
Swapan Kumar Ghosh of Ramakrishna Mission Vidyamandira for 
continuous inspiration and a number of valuable suggestions 
towards the improvement of the paper.

\end{document}